**Corner Regularity and Asymptotics for the Two-Dimensional Radiation Transport Equation**

Dean Wang

The Ohio State University

**Abstract**

The regularity of solutions to the two-dimensional radiation transport equation plays a central role in understanding the behavior of scalar flux near geometric corners. In a previous study of Kaper's corner singularity problem, the physical transport kernel was decomposed into its principal inverse-distance singularity and a lower-order remainder, reducing the corner asymptotic analysis to the question of whether the scalar flux possesses sufficient continuity to justify differentiation of the associated singular potential. The required regularity, however, remained unresolved. In this paper, we establish a sharp modulus-of-continuity estimate for the scalar flux generated by the physical transport kernel and establish the estimate $\omega_\phi = O(h\log(D/h))$, which implies Dini continuity. This estimate yields absolute convergence of the fluctuation term arising in the differentiated principal potential and permits a rigorous differentiation formula for transport solutions. Using this regularity framework, we derive the fixed-angle asymptotic behavior of the scalar flux and its first derivatives near a rectangular corner. For every fixed interior ray, the gradient exhibits logarithmic growth, while the scalar flux admits an $r\log(1/r)$ correction with an explicitly identified coefficient determined by the corner values of the flux and source. The analysis shows that a rectangular corner does not generate a singularity stronger than the logarithmic first-derivative singularity already present in smooth-boundary transport theory. Numerical calculations on benchmark problems are consistent with the theoretical predictions. These results provide a rigorous resolution of the regularity issue underlying Kaper's corner singularity problem and establish a mathematical framework for corner asymptotics in two-dimensional radiation transport.

## 1. Introduction

The regularity of solutions to the two-dimensional radiation transport equation is closely connected to the singular structure of the transport kernel. Even in domains with smooth boundaries, transport solutions generally fail to be smoothly differentiable up to the boundary. Previous work showed that first derivatives of the scalar flux may exhibit logarithmic growth, while higher derivatives possess progressively stronger boundary singularities governed by the distance to the boundary [1]. Thus, limited regularity is an intrinsic feature of the transport equation itself and not merely a consequence of geometric singularities.

The effect of corners on solution regularity is considerably less understood. In particular, Kaper's Problem 4 asks for the asymptotic behavior of the scalar flux and its derivatives near the corner of a quadrant for the one-group transport equation with isotropic scattering [3]. More broadly, one may ask whether a geometric corner produces singularities stronger than those already present in smooth-boundary transport theory. For elliptic equations, corner geometry can fundamentally alter local regularity and generate singular solutions that dominate the smooth-boundary behavior [4,5]. Whether an analogous phenomenon occurs for the physical two-dimensional radiation transport equation has remained an open question.

A major difficulty arises from the fact that the transport operator is itself singular. Consequently, geometric corner effects must be distinguished from the singularities already inherent in the transport kernel. In a recent study of Kaper's corner singularity problem [6], the physical transport kernel was decomposed into a principal inverse-distance singularity and a lower-order remainder. Differentiation of the associated transport potential separated the gradient into an explicit boundary contribution and an interior singular-integral term. The boundary contribution was shown to possess logarithmic growth near a rectangular corner, suggesting that the leading corner singularity might be no stronger than the logarithmic singularity already observed near smooth boundaries. However, control of the interior singular-integral term required the scalar flux to possess a Dini modulus of continuity. Establishing this regularity property was left as an open problem, preventing a rigorous characterization of the corner asymptotics.

The purpose of the present work is to resolve this outstanding regularity issue directly from the physical transport equation. We first establish a quantitative modulus-of-continuity estimate for the scalar flux,

$$\omega_\phi(h) = O\left(h\left(1 + \log\frac{D}{h}\right)\right), \tag{1.1}$$

which implies Dini continuity. The proof is obtained by decomposing the transport kernel into its principal inverse-distance component and a lower-order remainder and then deriving sharp continuity estimates for the associated integral operators. The resulting modulus is shown to be essentially sharp.

The Dini estimate plays a decisive role in the corner analysis. It guarantees absolute convergence of the fluctuation term appearing in the differentiated principal potential and provides a rigorous differentiation formula for transport solutions. This allows the singular contribution generated by the two sides meeting at a rectangular corner to be isolated explicitly. As a consequence, we derive the fixed-angle asymptotic behavior of both the scalar flux and its first derivatives near the corner.

For a rectangular corner $P$ and a fixed interior ray $x = re_\theta$, $0 < \theta < \pi/2$, we show that

$$\partial_{x_k}\phi(re_\theta) = A\log\frac{1}{r} + O_\theta(1),\ k = 1, 2, \tag{1.2}$$

where

$$A = \frac{\Sigma_s\phi(P) + Q(P)}{4}. \tag{1.3}$$

These results show that a rectangular corner does not generate a singularity stronger than the logarithmic first-derivative singularity already present in smooth-boundary transport theory. Instead, the corner determines the explicit coefficient of the leading singular term.

The significance of the present work extends beyond the particular corner geometry considered here. The Dini regularity theorem provides a mathematical framework for differentiating transport potentials with physical kernels and furnishes the regularity needed for corner singularity analysis in radiation transport. From a numerical perspective, the results provide an analytical basis for interpreting solution behavior near nonsmooth boundaries and for distinguishing genuine corner singularities from discretization artifacts.

The remainder of the paper is organized as follows. Section 2 formulates the transport integral equation and decomposes the physical kernel into its principal and remainder parts. Section 3 establishes continuity and Dini regularity of the scalar flux. Section 4 derives the corner asymptotics and identifies the leading logarithmic coefficient. Section 5 summarizes the main results. Section 6 presents numerical calculations that test the theoretical predictions, and Section 7 concludes the paper.

## 2. Integral Formulation and Kernel Decomposition

### 2.1. Integral equation

Let $G = (0, L_1) \times (0, L_2)$ be a rectangular domain and let $P = (0, 0)$ denote its lower-left corner. The scalar flux satisfies the integral equation [1]

$$\phi(x) = \int_G K\,(|x - y|)\phi(y)\,dy + F(x),\ x \in G. \tag{2.1}$$

Here $F(x) = \frac{1}{\Sigma_s}\int_G K(|x-y|)Q(y)dy$, where $Q$ is a bounded external source. $\Sigma_s > 0$ and $\Sigma_t > 0$ denote the scattering and total transport coefficients, respectively. It should be noted that vacuum boundary conditions are assumed throughout this study. As a result, the incoming boundary flux is identically zero and does not contribute to the regularity of the solution within the domain.

The physically relevant kernel is

$$K(q) = \frac{\Sigma_s}{2\pi q}\mathrm{Ki}_1(\Sigma_t q),\ q > 0, \tag{2.2}$$

where

$$\mathrm{Ki}_1(z) = \int_0^1 \frac{e^{-z/t}}{\sqrt{1-t^2}}\,dt \tag{2.3}$$

is the Bickley function of order one [1].

It is known that the integral equation admits a bounded solution [7]

$$\phi \in L^\infty(G). \tag{2.4}$$

For the corner asymptotics, we also assume that the external source is sufficiently smooth, for example that it belongs to $C^{1,\alpha}(\bar{G})$ for some $0 < \alpha < 1$. The source is assumed to be finite at the corner.

### 2.2. Decomposition of the transport kernel

To identify the singular structure of the transport operator, we first examine the behavior of the Bickley function for small arguments. Since the kernel singularity arises from the regime $|x-y| \to 0$, the asymptotic expansion of $\mathrm{Ki}_1(z)$ as $z \downarrow 0$ determines the leading behavior of the transport kernel.

**Lemma 2.1. Small-argument expansion of the Bickley function**

As $z \downarrow 0$, the small-argument expansion of the Bickley function, given in Eq. (2.3) has the form

$$\mathrm{Ki}_1(z) = \frac{\pi}{2} + z\log z + c_0 z + O\,(z^3), \tag{2.5}$$

where $c_0 = \gamma - 1 - \log 2$ and $\gamma$ denotes Euler's constant.

**Proof**

The Bickley functions satisfy the recurrence relation

$$\frac{d}{dz}\mathrm{Ki}_n(z) = -\mathrm{Ki}_{n-1}(z). \tag{2.6}$$

For $n = 1$,

$$\frac{d}{dz}\mathrm{Ki}_1(z) = -\mathrm{Ki}_0(z). \tag{2.7}$$

where $\mathrm{Ki}_0(z) = \int_0^1 \frac{e^{-z/t}}{t\sqrt{1-t^2}}\,dt$ is the modified Bessel function of the second kind, and its standard small-argument expansion gives [8]

$$\mathrm{Ki}_0(z) = -\log\left(\frac{z}{2}\right) - \gamma + O(z^2|\log z|), \quad z \downarrow 0. \tag{2.8}$$

Substituting Eq. (2.7) into (2.6) and integrating from 0 to $z$,

$$\mathrm{Ki}_1(z) - \mathrm{Ki}_1(0) = \int_0^z \left[\log\left(\frac{s}{2}\right) + \gamma + O(s^2|\log s|)\right] ds. \tag{2.9}$$

Since $\mathrm{Ki}_1(0) = \int_0^1 \frac{\mathrm{dt}}{\sqrt{1-t^2}} = \frac{\pi}{2}$, we obtain

$$\mathrm{Ki}_1(z) = \frac{\pi}{2} + \int_0^z \log s \, ds - (\log 2)z + \gamma z + O(z^3|\log z|). \tag{2.10}$$

Using $\int_0^z \log s \, ds = z \log z - z$, it follows that

$$\mathrm{Ki}_1(z) = \frac{\pi}{2} + z\log z + (\gamma - 1 - \log 2)z + O(z^3|\log z|). \tag{2.11}$$

This proves the stated expansion. □

Substituting $z = \Sigma_t q$ into Lemma 2.1 and using the kernel definition (2.2), we obtain

$$K(q) = \frac{\Sigma_s}{2\pi q}\left[\frac{\pi}{2} + \Sigma_t q \log(\Sigma_t q) + c_0 \Sigma_t q + O(q^3|\log q|)\right]. \tag{2.12}$$

Hence,

$$K(q) = \frac{\Sigma_s}{4q} + O\left(1 + |\log q|\right) \tag{2.13}$$

as $q \downarrow 0$.

The inverse-distance term is therefore the unique nonintegrable part of the kernel derivative and determines the leading singular behavior of the transport operator.

Motivated by this expansion, we decompose the kernel as

$$K(q) = K_0(q) + R(q), \tag{2.14}$$

where

$$K_0(q) = \frac{\Sigma_s}{4q}, \tag{2.15}$$

and

$$R(q) = K(q) - K_0(q). \tag{2.16}$$

The corresponding integral operators are

$$T_0 u(x) = \frac{\Sigma_s}{4}\int_G \frac{u(y)}{|x-y|} \, dy, \tag{2.17}$$

and

$$T_R u(x) = \int_G R\left(|x-y|\right) u(y) \, dy. \tag{2.18}$$

Writing $T_K = T_0 + T_R$, the integral equation (2.1) becomes

$$\phi = T_0\phi + T_R\phi + F. \tag{2.19}$$

Using Lemma 2.1 once more gives

$$R(q) = \frac{\Sigma_s \Sigma_t}{2\pi}\left[\log(\Sigma_t q) + c_0\right] + O\left(q^2|\log q|\right) \tag{2.20}$$

as $q \downarrow 0$. In particular,

$$|R(q)| \leq C(1 + |\log q|). \tag{2.21}$$

To estimate ￼, we differentiate the exact kernel rather than the asymptotic remainder. From (2.2),

$$K'(q) = -\frac{\Sigma_s}{2\pi q^2}\mathrm{Ki}_1(\Sigma_t q) + \frac{\Sigma_s \Sigma_t}{2\pi q}\mathrm{Ki}_1'(\Sigma_t q). \tag{2.22}$$

Using the small-argument expansions (2.8) and (2.11), together with $\mathrm{Ki}_1'(z) = -\mathrm{Ki}_0(z)$ gives

$$K'(q) = \frac{\Sigma_s}{4q^2} + \frac{\Sigma_s \Sigma_t}{2\pi q} + O(q|\log q|). \tag{2.23}$$

Since $K_0'(q) = -\frac{\Sigma_s}{4q^2}$, the leading $q^{-2}$ terms cancel in $R' = K' - K_0'$, yielding

$$R'(q) = \frac{\Sigma_s \Sigma_t}{2\pi q} + O(q|\log q|). \tag{2.24}$$

Therefore, for sufficiently small $q > 0$,

$$|R'(q)| \leq \frac{C}{q}. \tag{2.25}$$

Because $G$ is bounded, the constants may be enlarged so that these estimates remain valid throughout the range of distances occurring in the domain.

The logarithmic singularity of $R$ is locally integrable in two dimensions, and the derivative estimate is also integrable with respect to area near the origin. These observations form the basis for the regularity results established in the next section.

## 3. Regularity of the Scalar Flux

### 3.1. Continuity estimate for the principal operator

In this section, $u$ denotes a generic bounded density, while $\phi$ denotes the scalar-flux solution of the transport equation. The operator estimates are first proved for $u$ and are subsequently applied to $u = \phi$.

For a bounded function $u$, define its modulus of continuity, when $u$ is continuous, by

$$\omega_u(h) = \sup_{\substack{x,x' \in \bar{G} \\ |x - x'| \leq h}} |u(x) - u(x')|. \tag{3.1}$$

We first estimate the difference of two values of the principal potential.

**Lemma 3.1.** There is a constant $C > 0$, depending only on the domain and the kernel coefficient, such that for $u \in L^\infty(G)$ and sufficiently small $h = |x - x'| > 0$,

$$|T_0 u(x) - T_0 u(x')| \leq C\|u\|_{L^\infty(G)} h \left(1 + \log\frac{D}{h}\right), \tag{3.2}$$

where $D$ is any fixed length scale larger than the diameter of $G$.

**Proof**

By the definition of $T_0$,

$$|T_0 u(x) - T_0 u(x')| \leq \frac{\Sigma_s}{4}\|u\|_{L^\infty(G)} \int_G \left|\frac{1}{|x - y|} - \frac{1}{|x' - y|}\right| dy. \tag{3.3}$$

Split the integration domain into the near region

$$N = \{y \in G: |x-y| \le 2h\}$$

and its complement $G\backslash N$.

For $y \in N$,

$$|x'-y| \le |x'-x| + |x-y| \le 3h. \tag{3.4}$$

Hence

$$\int_N \left|\frac{1}{|x-y|} - \frac{1}{|x'-y|}\right| dy \le \int_N \left(\frac{1}{|x-y|} + \frac{1}{|x'-y|}\right) dy \le Ch. \tag{3.5}$$

For $y \in G\backslash N$, $|x-y| > 2h$. Applying the mean value theorem to the function

$$f(z) = \frac{1}{|z-y|}, \tag{3.6}$$

and using

$$|\nabla f(z)| = \tfrac{1}{|z-y|^2}, \tag{3.7}$$

gives

$$\left|\frac{1}{|x-y|} - \frac{1}{|x'-y|}\right| \le \frac{Ch}{|x-y|^2}. \tag{3.8}$$

Therefore,

$$\int_{G\backslash N} \left|\frac{1}{|x-y|} - \frac{1}{|x'-y|}\right| dy \le Ch \int_{G\backslash N} \frac{dy}{|x-y|^2} \le Ch \int_{2h}^{D} \frac{d\rho}{\rho} \le Ch \log\frac{D}{h}\,. \tag{3.9}$$

Combining (3.5) and (3.9) proves (3.2). □

The logarithmic factor arises from the annular integral of the derivative of the inverse-distance kernel. It is characteristic of this two-dimensional singularity.

**Remark 3.1.** The logarithmic factor in (3.2) is optimal in general, even for a constant density. Let $u \equiv 1$, and let $x_0 = (0, a)$, $0 < a < L_2$, lie in the relative interior of the side $x_1 = 0$, and set $x_h = (h, a)$, $h > 0$. Direct integration over the rectangle gives

$$T_0 1(x_h) \;= \frac{\Sigma_s}{4} \int_0^{L_1} \left[\operatorname{arsinh}\frac{a}{|y_1-h|} + \operatorname{arsinh}\frac{L_2-a}{|y_1-h|}\right] dy_1\,. \tag{3.10}$$

Since

$$\operatorname{arsinh}\frac{a}{s} + \operatorname{arsinh}\frac{L_2-a}{s} = 2\log\frac{1}{s} + O(1), \qquad s \downarrow 0, \tag{3.11}$$

comparison with the corresponding expression at $x_0$ yields

$$T_0 1(x_h) - T_0 1(x_0) = \frac{\Sigma_s}{2} h \log\frac{1}{h} + O(h)\,, \qquad h \downarrow 0. \tag{3.12}$$

Thus the modulus $h(1 + \log(D/h))$ in Lemma 3.1 cannot, in general, be replaced by a smaller modulus uniformly over bounded density. The logarithmic loss from Lipschitz continuity is therefore intrinsic to the boundary behavior of the inverse-distance potential and occurs even for smooth densities.

### 3.2. Lipschitz estimate for the exact remainder

**Lemma 3.2.**
Let

$$(T_R u)(x) = \int_G R(|x-y|)u(y)\,dy\,, \tag{3.13}$$

where $R$ is the remainder kernel defined by (2.16). Then there exists a constant $C > 0$ such that

$$|(T_R u)(x) - (T_R u)(x')| \le C\|u\|_{L^\infty(G)}|x - x'|, \quad x, x' \in \bar{G}. \tag{3.14}$$

Consequently, $T_R: L^\infty(G) \to C^{0,1}(\bar{G})$.

**Proof.**
By (2.20),

$$|R'(q)| \le \frac{C}{q}, \quad 0 < q \le D, \tag{3.15}$$

where $D > \operatorname{diam}(G)$.
Hence

$$|\nabla_x R(|x-y|)| = \left|R'(|x-y|)\frac{x-y}{|x-y|}\right| \le \frac{C}{|x-y|}. \tag{3.16}$$

let $h = |x - x'|$, where $x_t = (1-t)x + tx'$, $0 \le t \le 1$.
For almost every $y \in G$, the fundamental theorem of calculus gives

$$|R(|x-y|) - R(|x'-y|)| = \int_0^1 \nabla_x R(|x_t - y|) \cdot (x - x')\,dt\,. \tag{3.17}$$

Using (3.16).

$$|R(|x-y|) - R(|x'-y|)| \le Ch \int_0^1 \frac{dt}{|x_t - y|}. \tag{3.18}$$

Therefore,

$$|(T_R u)(x) - (T_R u)(x')| \le \|u\|_{L^\infty(G)} \int_G |R(|x-y|) - R(|x'-y|)|\,dy$$

$$\le Ch\|u\|_{L^\infty(G)} \int_G \int_0^1 \frac{dt}{|x_t - y|}\,dy\,. \tag{3.19}$$

Applying Fubini's theorem,

$$|(T_R u)(x) - (T_R u)(x')| \le Ch\|u\|_{L^\infty(G)} \int_0^1 \left(\int_G \frac{dy}{|x_t - y|}\right) dt\,. \tag{3.20}$$

Since $G \subset B(x_t, D)$,

$$\int_G \frac{dy}{|x_t - y|} \le \int_{B(x_t,D)} \frac{dy}{|x_t - y|} = 2\pi \int_0^D d\rho = 2\pi D, \tag{3.21}$$

uniformly in $t$. Substituting (3.18) into (3.17) yields

$$|(T_R u)(x) - (T_R u)(x')| \le Ch\|u\|_{L^\infty(G)} = C\|u\|_{L^\infty(G)}|x - x'|, \tag{3.22}$$

Thus $T_R u$ is Lipschitz continuous on $\bar{G}$. □

**Remark 3.2.** The Lipschitz exponent in Lemma 3.2 is optimal in general in the Hölder scale. Indeed, whenever $T_R u$ is differentiable at a point $x_0$ in a direction $e$ and $e \cdot \nabla T_R u(x_0) \neq 0$, one has

$$T_R u(x_0 + he) - T_R u(x_0) = he \cdot \nabla T_R u(x_0) + o(h), \qquad h \downarrow 0. \tag{3.23}$$

so no uniform modulus $o(h)$ can hold on a class of densities containing such an input. Thus the exponent 1 in Lemma 3.2 cannot, in general, be improved. In contrast with the principal operator, however, the remainder produces no logarithmic loss: its weaker kernel singularity yields a Lipschitz potential rather than the $h\log(1/h)$ boundary modulus associated with $T_0$.

### 3.3. Dini continuity

A modulus of continuity $\omega$ is called a Dini modulus [9,10] if

$$\int_0^1 \frac{\omega(h)}{h}\, dh < \infty. \tag{3.24}$$

The logarithmic modulus

$$\omega_*(h) = h\left(1 + \log\frac{D}{h}\right) \tag{3.25}$$

satisfies the Dini condition. Indeed,

$$\int_0^\delta \frac{\omega_*(h)}{h}\, dh = \int_0^\delta \left(1 + \log\frac{D}{h}\right) dh < \infty. \tag{3.26}$$

Dini moduli play a fundamental role in singular-integral regularity theory [10].

**Theorem 3.1.**

Suppose that $\phi \in L^\infty(G)$ solves (2.19), and that $F$ is generated by a bounded external source $Q$, as given in Section 2.1. Then

$$|\phi(x) - \phi(x')| \leq C\|\phi\|_{L^\infty(G)} h\left(1 + \log\frac{D}{h}\right) + \omega_F(h),\ h = |x - x'|, \tag{3.27}$$

where

$$\omega_F(h) = \sup_{\substack{x,x' \in \bar{G} \\ |x - x'| \leq h}} |F(x) - F(x')|. \tag{3.28}$$

In particular, the source-potential structure gives the same logarithmic modulus for $F$, so $\phi$ is Dini continuous; no Lipschitz assumption on $F$ is needed.

**Proof**

From Eq. (2.19),

$$\phi = T_0\phi + T_R\phi + F. \tag{3.29}$$

Applying Lemmas 3.1 and 3.2 yields

$$|\phi(x) - \phi(x')| \leq |(T_0\phi)(x) - (T_0\phi)(x')| + |(T_R\phi)(x) - (T_R\phi)(x')| + |F(x) - F(x')|$$

$$\leq C\|\phi\|_{L^\infty(G)} h\left(1+\log\frac{D}{h}\right) + C\|\phi\|_{L^\infty(G)} h + \omega_F(h). \tag{3.30}$$

Since

$$h \leq h\left(1+\log\frac{D}{h}\right), \quad 0 < h < D, \tag{3.31}$$

the Lipschitz term may be absorbed into the logarithmic modulus by enlarging the constant $C$. This gives (3.27).

To estimate the source term, decompose $K = K_0 + R$ and use the same estimates as in Lemmas 3.1 and 3.2. Since $F = (1/\Sigma_s)(T_0 Q + T_R Q)$,

$$\omega_F(h) \leq C\|Q\|_{L^\infty(G)} h\left(1+\log\frac{D}{h}\right), \quad 0 < h \ll 1. \tag{3.32}$$

Combining this with (3.27) and the bound for the remainder operator gives

$$\omega_\phi(h) \leq C\left(\|\phi\|_{L^\infty(G)} + \|Q\|_{L^\infty(G)}\right) h\left(1+\log\frac{D}{h}\right). \tag{3.33}$$

The modulus $h\left(1+\log\frac{D}{h}\right)$ is Dini because $\int_0^\delta \left(1+\log\frac{D}{h}\right) dh$ is finite. Thus $\phi$ has a Dini-continuous representative on the closure of $G$. □

## 4. Corner Asymptotics

### 4.1. Differentiation formula for the principal potential

The leading corner behavior is determined by the principal potential

$$T_0 u(x) = \frac{\Sigma_s}{4} \int_G \frac{u(y)}{|x-y|}\, dy. \tag{4.1}$$

Because the gradient of the kernel behaves like $|x-y|^{-2}$, direct differentiation is not absolutely integrable. The appropriate differentiation formula is obtained by separating the singular part from the fluctuation of the density.

**Lemma 4.1.**

Let $u \in C(\bar{G})$ and suppose its modulus of continuity satisfies a Dini condition, $\int_0^1 \frac{\omega_u(\rho)}{\rho} d\rho < \infty$. Then $T_0 u \in C^1(G)$, and for every $x \in G$,

$$\partial_{x_k} T_0 u(x) = -\frac{\Sigma_s}{4} \int_G \frac{x_k - y_k}{|x-y|^3} [u(y) - u(x)]\, dy - \frac{\Sigma_s}{4} u(x) \int_{\partial G} \frac{n_k(y)}{|x-y|}\, dS_y,\ k = 1,2, \tag{4.2}$$

where $n(y)$ is the outward unit normal to the boundary. The integral containing $u(y) - u(x)$ is absolutely convergent.

**Proof**

The distributional derivative of the kernel $|z|^{-1}$ in two dimensions is

$$\partial_{z_k} |z|^{-1} = -\,\mathrm{p.v.}\frac{z_k}{|z|^3}, \tag{4.3}$$

where p. v. denotes the Cauchy principal value. Since the average of the $k$-th component of the unit normal over a circle vanishes,

$$\int_{|z|=\varepsilon} \frac{z_k}{|z|} dS(z) = 0, \tag{4.4}$$

no point-mass term appears.

Therefore,

$$\partial_{x_k} T_0 u(x) = -\frac{\Sigma_s}{4} \text{p.v.} \int_G \frac{x_k - y_k}{|x-y|^3} \; u(y)\, dy\,. \tag{4.5}$$

To remove the principal value singularity, write

$$u(y) = u(x) + \big(u(y) - u(x)\big). \tag{4.6}$$

Substituting into (4.3) yields

$$\partial_{x_k} T_0 u(x) \;\; = -\frac{\Sigma_s}{4} \int_G \frac{x_k - y_k}{|x-y|^3} \big(u(y) - u(x)\big)\, dy - \frac{\Sigma_s}{4} u(x) \;\text{p.v.} \int_G \frac{x_k - y_k}{|x-y|^3} dy\,. \tag{4.7}$$

The first integral is absolutely convergent. Indeed,

$$\left| \frac{x_k - y_k}{|x-y|^3} [u(y) - u(x)] \right| \le \frac{\omega_u(|x-y|)}{|x-y|^2}, \tag{4.8}$$

and hence

$$\int_{|x-y|<\delta} \frac{\omega_u(|x-y|)}{|x-y|^2} dy = 2\pi \int_0^\delta \frac{\omega_u(\rho)}{\rho} d\rho\,. \tag{4.9}$$

The right-hand side is finite by the Dini condition.

It remains to evaluate the constant-density term. Let $G_\varepsilon = G \backslash \overline{B_\varepsilon(x)}$.

Since

$$\nabla_y |x-y|^{-1} = \frac{x-y}{|x-y|^3}, \tag{4.10}$$

the divergence theorem gives

$$\int_{G_\varepsilon} \frac{x_k - y_k}{|x-y|^3} dy = \int_{\partial G} \frac{n_k(y)}{|x-y|} dS(y) + \int_{\partial G_\varepsilon} \frac{\nu_k(y)}{|x-y|} dS(y)\,, \tag{4.11}$$

where $\nu$ is the outward normal of $G_\varepsilon$.

On the inner boundary, $\nu(y) = -\frac{y-x}{|y-x|}$, and

$$\int_{\partial G_\varepsilon} \frac{\nu_k(y)}{|x-y|} dS(y) = -\int_0^{2\pi} (\cos\theta\,, \sin\theta)_k d\theta = 0. \tag{4.12}$$

Letting $\varepsilon \downarrow 0$ therefore yields

$$\text{p.v.} \int_G \frac{x_k - y_k}{|x-y|^3} \, dy = \int_{\partial G} \frac{n_k(y)}{|x-y|} \;\; dS(y)\,. \tag{4.13}$$

Substituting (4.13) into (4.7) gives (4.2).

Finally, continuity of the right-hand side follows from the Dini estimate (4.5) together with dominated convergence away from the singularity. Hence $T_0 u \in C^1(G)$. □

**4.2. The boundary logarithm**

Let

$$x = re_\theta,\ e_\theta = (\cos\theta, \sin\theta),\ 0 < \theta < \frac{\pi}{2}, \tag{4.14}$$

where $\theta$ is fixed and $r \downarrow 0$.

The boundary term in (4.2) contains the two sides meeting at the corner. We first examine the boundary integral. Consider $k = 1$. The only side producing an unbounded contribution as $r \downarrow 0$ is the vertical side $y_1 = 0$, whose outward normal satisfies $n_1 = -1$. Its contribution is

$$-\int_0^{L_2} \frac{ds}{\sqrt{r^2\cos^2\theta + (s - r\sin\theta)^2}}. \tag{4.15}$$

Let $a = r\cos\theta$, $b = r\sin\theta$. Using

$$\int \frac{ds}{\sqrt{a^2 + (s-b)^2}} = \operatorname{arsinh}\left(\frac{s-b}{a}\right), \tag{4.16}$$

we obtain

$$\int_0^{L_2} \frac{ds}{\sqrt{a^2 + (s-b)^2}} = \operatorname{arsinh}\left(\frac{L_2 - b}{a}\right) + \operatorname{arsinh}\left(\frac{b}{a}\right). \tag{4.17}$$

Since $a = r\cos\theta$,

$$\operatorname{arsinh}\left(\frac{L_2 - b}{a}\right) = \log\frac{2L_2}{r\cos\theta} + O(r), \tag{4.18}$$

while

$$\operatorname{arsinh}\left(\frac{b}{a}\right) = \operatorname{arsinh}(\tan\theta) \tag{4.19}$$

is independent of $r$. Therefore

$$\int_0^{L_2} \frac{ds}{\sqrt{r^2\cos^2\theta + (s - r\sin\theta)^2}} = \log\frac{1}{r} + O_\theta(1), \tag{4.20}$$

where $O_\theta(1)$ denotes a quantity bounded as $r \to 0$, where the bound may depend on the fixed angle $\theta$.

Hence

$$\int_{\partial G} \frac{n_1(y)}{|re_\theta - y|}\, dS_y = -\log\frac{1}{r} + O_\theta(1). \tag{4.21}$$

An identical calculation for the horizontal side $y_2 = 0$ gives

$$\int_{\partial G} \frac{n_2(y)}{|re_\theta - y|}\, dS_y = -\log\frac{1}{r} + O_\theta(1). \tag{4.22}$$

The remaining sides stay a positive distance away from the corner and therefore contribute only bounded terms. Next consider the fluctuation integral in (4.2). We now specialize Lemma 4.1 to the scalar-flux solution $\phi$. By Theorem 3.1, $\phi$ is Dini continuous on $\bar{G}$, and therefore Lemma 4.1 applies with $u = \phi$. Thus,

$$\omega_\phi(\rho) \leq C\rho\left(1 + \log\frac{D}{\rho}\right). \tag{4.23}$$

Consequently,

$$\left|\iint_{|x-y|<\delta} \frac{x_k - y_k}{|x-y|^3}[\phi(y) - \phi(x)]\,dy\right| \leq \int_{|x-y|<\delta} \frac{|\phi(y) - \phi(x)|}{|x-y|^2}\,dy$$
$$\leq C\int_0^\delta \frac{\omega_\phi(\rho)}{\rho}\,d\rho < \infty, \tag{4.24}$$

and the far-field contribution is also bounded. Thus

$$\int_G \frac{x_k - y_k}{|x-y|^3}[\phi(y) - \phi(x)]\,dy = O_\theta(1). \tag{4.25}$$

Combining the boundary estimates (4.21) and (4.22) with the bounded fluctuation estimate (4.25) in the differentiation formula (4.2), we obtain

$$\partial_{x_k} T_0\phi(re_\theta) = \frac{\Sigma_s}{4}\phi(re_\theta)\log\frac{1}{r} + O_\theta(1), \qquad k = 1,2. \tag{4.26}$$

Since $\phi$ is continuous at the corner,

$$\phi(re_\theta) = \phi(P) + O\left(\omega_\phi(r)\right), \tag{4.27}$$

hence

$$\phi(re_\theta)\log\frac{1}{r} = \phi(P)\log\frac{1}{r} + O\left(\omega_\phi(r)\log\frac{1}{r}\right). \tag{4.28}$$

Using

$$\omega_\phi(r) = O\left(r\log\frac{D}{r}\right), \tag{4.29}$$

gives

$$\omega_\phi(r)\log\frac{1}{r} = O\left(r\log^2\frac{1}{r}\right) \to 0. \tag{4.30}$$

Therefore

$$\phi(re_\theta)\log\frac{1}{r} = \phi(P)\log\frac{1}{r} + o(1). \tag{4.31}$$

Hence the coefficient in (4.26) may be replaced by $\phi(P)$, giving

$$\partial_{x_k} T_0\phi(re_\theta) = \frac{\Sigma_s}{4}\phi(P)\log\frac{1}{r} + O_\theta(1), \quad k = 1,2. \tag{4.32}$$

It remains to estimate the contribution of the remainder operator.

By (3.16), i.e., $|\nabla_x R(|x-y|)| \leq \frac{C}{|x-y|}$. For $k = 1, 2$,

$$\left|\partial_{x_k}(R|x-y|)\phi(y)\right| \le C\|\phi\|_{L^\infty(G)}\frac{1}{|x-y|}. \tag{4.33}$$

Since $|x-y|^{-1} \in L^1(G)$ in two dimensions, the right-hand side is integrable over $G$. Therefore, the dominated convergence theorem permits differentiation under the integral sign, yielding

$$\partial_{x_k}T_R\phi(x) = \int_G \partial_{x_k}R(x-y)\phi(y)dy\,, \quad k = 1, 2. \tag{4.34}$$

Consequently,

$$|\nabla T_R\phi(x)| \le \|\phi\|_{L^\infty(G)}\int_G |\nabla_x R(x-y)|dy \le C\|\phi\|_{L^\infty(G)}\int_G \frac{dy}{|x-y|}. \tag{4.35}$$

Because $G$ is bounded,

$$\int_G \frac{dy}{|x-y|} \le \int_{B(x,D)} \frac{dy}{|x-y|} = 2\pi\int_0^D d\rho = 2\pi D, \tag{4.36}$$

Where $D > \mathrm{diam}(G)$. Hence

$$|\nabla T_R\phi(x)| = O(1), \tag{4.37}$$

and, in particular

$$|\nabla T_R\phi(re_\theta)| = O(1), \tag{4.38}$$

Thus the exact remainder does not contribute to the leading logarithmic term.

**4.3. The source contribution**

Because $F(x) = \frac{1}{\Sigma_s}\int_G K(|x-y|)Q(y)dy$, its principal part is $\frac{1}{4}\int_G \frac{Q(y)}{|x-y|}dy$. Since $Q \in C^{1,\alpha}(G)$ by assumption, the same differentiation formula applies to $Q$, and the remainder potential has bounded first derivatives. Consequently, for each fixed interior angle $\theta$, the source contribution to each coordinate derivative has leading coefficient $Q(P)/4$.

It follows that, for a sufficiently smooth source,

$$\partial_{x_k}F(re_\theta) = \frac{Q(P)}{4}\log\frac{1}{r} + O_\theta(1),\ k = 1, 2. \tag{4.39}$$

Here $Q(P)$ denotes the external-source value at the corner $P$.

Set $A = [\Sigma_s\phi(P) + Q(P)]/4$. Combining the scalar-flux and source contributions gives

$$\partial_{x_1}\phi(re_\theta) = A\log\frac{1}{r} + O_\theta(1), \tag{4.40}$$

and

$$\partial_{x_2}\phi(re_\theta) = A\log\frac{1}{r} + O_\theta(1). \tag{4.41}$$

These formulas identify the leading corner singularity of the gradient of the scalar flux.

**4.4. Radial asymptotics**

Along the ray $x = re_\theta$, the chain rule gives

$$\frac{d}{dr}\phi(re_\theta) = \cos\theta\ \partial_{x_1}\phi(re_\theta) + \sin\theta\ \partial_{x_2}\phi(re_\theta). \tag{4.42}$$

Using (4.40) and (4.41), we obtain the radial derivative

$$\frac{d}{dr}\phi(re_\theta) = A(\cos\theta + \sin\theta)\log\frac{1}{r} + O_\theta(1). \tag{4.43}$$

Since $\phi$ is continuous at $P$, integration from 0 to $r$ gives

$$\phi(re_\theta) - \phi(P) = A(\cos\theta + \sin\theta)\int_0^r \log\frac{1}{s}\, ds + O_\theta(r). \tag{4.44}$$

The elementary integral satisfies

$$\int_0^r \log\frac{1}{s}\, ds = r\log\frac{1}{r} + r. \tag{4.45}$$

The term proportional to $r$ is absorbed into the $O_\theta(r)$ remainder. We therefore obtain the fixed-angle corner asymptotic

$$\phi(re_\theta) - \phi(P) = A(\cos\theta + \sin\theta) r\log\tfrac{1}{r} + O_\theta(r). \tag{4.46}$$

The asymptotic is asserted for each fixed angle strictly inside the corner. The bound on the remainder is not claimed to be uniform as the angle approaches either side of the rectangle.

## 5. Main Result

### Theorem 5.1 (Corner Asymptotics and Regularity)

Let $G = (0, L_1) \times (0, L_2)$, and suppose that the scalar flux $\phi \in L^\infty(G)$ solves the two-dimensional radiation transport integral equation (2.1) with the physical kernel (2.2). Assume that the external source is sufficiently smooth, for example $C^{1,\alpha}(\bar{G})$ for some $0 < \alpha < 1$, and that the source and solution satisfy the boundedness conditions used above. Then the following conclusions hold.

1. The scalar flux has a continuous representative on $\bar{G}$ with logarithmic modulus of continuity

   $$|\phi(x) - \phi(x')| \leq Ch\left(1 + \log\frac{D}{h}\right),\ h = |x - x'|, \tag{5.1}$$

   for sufficiently small $h$. In particular, $\phi$ is Dini continuous.

2. For every fixed interior angle $0 < \theta < \pi/2$, the coordinate derivatives satisfy

   $$\partial_{x_k}\phi(re_\theta) = A\log\frac{1}{r} + O_\theta(1),\ k = 1, 2, \tag{5.2}$$

   where $A$ is one quarter of the sum of the scattering coefficient multiplied by $\phi(P)$ and the external-source value at $P$.

3. The scalar flux has the radial expansion

   $$\phi(re_\theta) - \phi(P) = A(\cos\theta + \sin\theta) r\log\frac{1}{r} + O_\theta(r). \tag{5.3}$$

**Proof**

The continuity estimate and Dini regularity follow from Theorem 3.1. Lemma 4.1 supplies the differentiation formula for the principal potential, with absolute convergence of the fluctuation integral ensured by Dini continuity. The boundary estimates (4.21) – (4.22) generate the logarithmic contribution, while the fluctuation term and the derivative of the exact remainder remain bounded. The source potential contributes its own leading logarithmic term. This proves (5.2). Integrating the directional derivative along a fixed ray yield (5.3). □

**Remark.** The kernel decomposition $K = K_0 + R$, shows that the inverse-distance principal kernel is solely responsible for the leading logarithmic behavior, while the exact remainder contributes only bounded corrections to the gradient. Consequently, the full transport kernel possesses the same leading fixed-angle corner asymptotics as its principal part. The source potential contributes a logarithmic term of the same order as the scattering contribution.

An important consequence of Theorem 5.1 is that the corner exhibits the same leading logarithmic first-derivative singularity already known for smooth-boundary transport theory. The principal contribution of this work is the derivation of a quantitative continuity estimate and a fixed-angle corner asymptotic expansion directly from the physical transport kernel, thereby providing a rigorous framework for understanding corner effects in two-dimensional radiation transport problems and their impact on numerical approximation near nonsmooth boundaries.

In addition, the leading logarithmic behavior in (5.2) formally suggests $D^m\phi(re_\theta) = O_\theta(r^{1-m}), m \geq 2$. This higher-order behavior cannot, however, be obtained rigorously by direct differentiation of (5.2), since no derivative bounds are available for its $O_\theta(1)$ remainder. Establishing higher-order asymptotics would require differentiating the underlying potential representation and controlling the resulting fluctuation and remainder terms. Such an analysis is beyond the scope of the present work.

## 6. Numerical Results

We examine the corner asymptotics using the unit-square benchmark on $G = (0,5)^2$, with reflective boundary conditions on the left and bottom sides and vacuum boundary conditions on the top and right sides. A uniform source, $Q = 1$, is prescribed throughout the domain. The relatively large domain is chosen to minimize interference between boundary effects from different sides. To reduce computational cost, reflective boundary conditions are imposed on the left and bottom boundaries, exploiting symmetry to restrict the computation to one quadrant of an equivalent larger problem with vacuum boundaries. Since the asymptotic analysis is local and concerns the behavior of the solution in the immediate vicinity of the corner, these reflective boundaries are not expected to affect the leading-order corner coefficient.

Two cases are considered: $(\Sigma_t, \Sigma_s) = (1, 0)$, corresponding to pure absorption, and $(1, 0.8)$, which includes scattering. The angular variable is discretized using the $S_{20}$ level-symmetric quadrature, while the spatial discretization employs the diamond-difference (DD) scheme on a fine $3000 \times 3000$ mesh. Figure 1 shows the resulting scalar-flux distributions for the two cases.

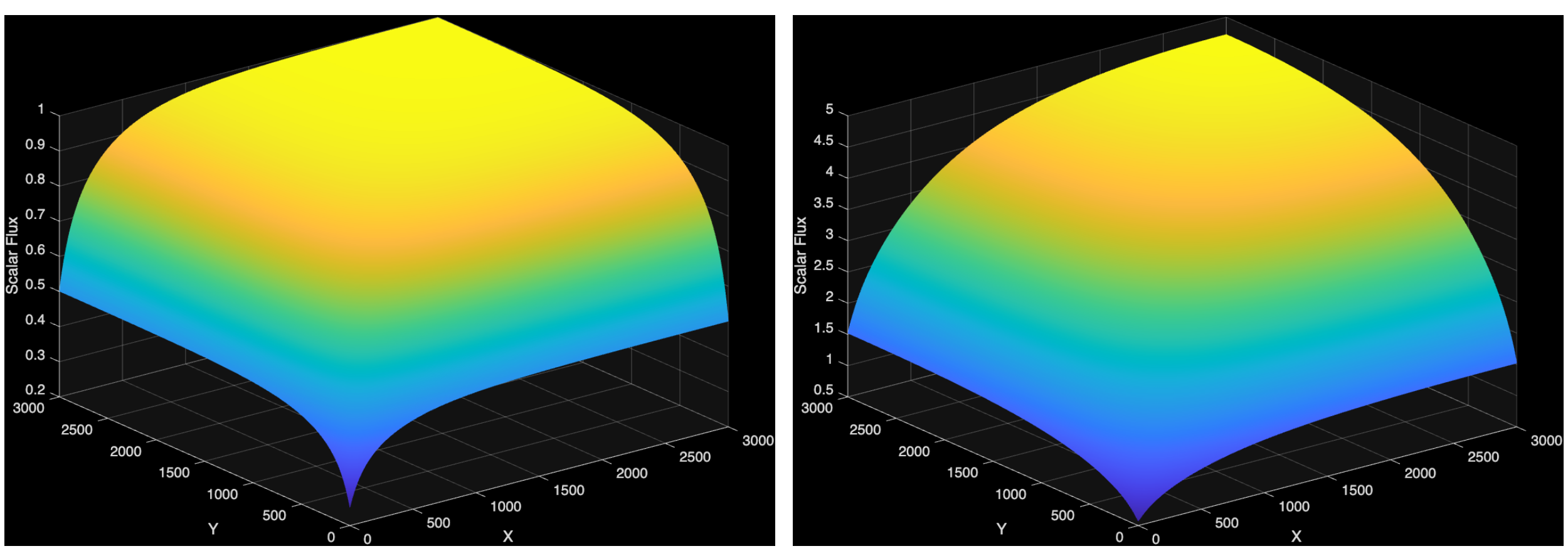

(a) (b)

Figure 1. Scalar-flux distributions for (a) Case 1 and (b) Case 2.

Along the diagonal ray $x = re_\theta$, the asymptotic expansion in Eq. (5.3) takes the form

$$\phi(re_\theta) - \phi(P) = A(\cos\theta + \sin\theta) r \log\frac{1}{r} + O_\theta(r),$$

where $A$ is the predicted leading coefficient, $A = \frac{\Sigma_s \phi(P) + Q(P)}{4}$, and $P$ is the corner.

For Case 1, $\Sigma_s = 0$ and $Q(P) = 1$, so the predicted coefficient is $A_{pred} = 0.250$. For Case 2, the computed corner value $\phi(P) = 0.56871596$ gives $A_{pred} \approx 0.364$.

To examine the asymptotic prediction at finite distances from the corner, we define the numerical estimator

$$A_{num} = \frac{\phi(re_\theta) - \phi(P)}{(\cos\theta + \sin\theta) r \log\frac{1}{r}}. \tag{6.1}$$

Thus, the finite-distance estimator is expected to approach the analytical coefficient as the corner is approached. The convergence can nevertheless be slow because the normalized contribution of the $O_\theta(r)$ remainder decays only as $1/\log(1/r)$.

Tables 1 and 2 compare the predicted coefficient with the numerical estimates for two interior rays: $\theta = \pi/4$ and $\theta = \text{atan}\frac{1}{2} \approx 0.1476\pi$, respectively.

**Table 1.** Analytical and numerical estimates of $A$ along the ray $\theta = \pi/4$.

| **Distance to corner *P*** $r$ | **Case 1** | | **Case 2** | |
|---|---|---|---|---|
| | $A_{pred}$ | $A_{num}$ | $A_{pred}$ | $A_{num}$ |
| 0.00942809 | 0.250 | 0.224 | 0.364 | 0.387 |
| 0.01414214 | | 0.244 | | 0.424 |
| 0.01885618 | | 0.259 | | 0.455 |
| 0.02357023 | | 0.272 | | 0.481 |
| 0.02828427 | | 0.284 | | 0.505 |
| 0.03299832 | | 0.294 | | 0.527 |

**Table 2.** Analytical and numerical estimates of $A$ along the ray $\theta = \text{atan}\frac{1}{2}$.

| **Distance to corner *P*** $r$ | **Case 1** | | **Case 2** | |
|---|---|---|---|---|
| | $A_{pred}$ | $A_{num}$ | $A_{pred}$ | $A_{num}$ |
| 0.00931695 | 0.250 | 0.219 | 0.364 | 0.379 |
| 0.01304373 | | 0.234 | | 0.408 |
| 0.01677051 | | 0.246 | | 0.432 |
| 0.02049729 | | 0.257 | | 0.454 |
| 0.02422407 | | 0.267 | | 0.474 |
| 0.02795085 | | 0.276 | | 0.492 |

In Case 1, the estimates at the smallest distances, 0.224 and 0.219, lie below the predicted coefficient 0.250, suggesting that the angular discretization error may still be significant in this nonscattering case and that a higher-order quadrature set could improve the agreement. In Case 2, the corresponding estimates, 0.387 and 0.379, lie slightly above the predicted coefficient 0.364, indicating that further spatial refinement may be beneficial. Nevertheless, these discrepancies are consistent with finite-distance corrections. It is noted that increasing the order of the $S_N$ quadrature improves the agreement with the theoretical predictions. For both cases, the numerical estimates increase as $r$ increases. This trend is consistent with the presence of a subleading $O_\theta(r)$ term in Eq. (5.3). After normalization by $1/\log(1/r)$, this term contributes $O_\theta(1/\log(1/r))$, which vanishes only logarithmically as $r \to 0$.

It would be of interest to investigate whether the leading coefficient can be recovered more accurately by further refining the mesh and increasing the angular resolution. Since $S_{20}$ is the highest-order level-symmetric quadrature available, this would require a different quadrature set. This topic will be explored in future work.

**7. Conclusions**

The principal new result is the proof of Dini continuity for the scalar flux generated by the physical transport kernel. This regularity closes the outstanding gap in the earlier corner-singularity analysis and permits rigorous derivation of the fixed-angle corner asymptotics. By decomposing the physical transport kernel into a principal inverse-distance singularity and an exact remainder, we showed that the principal operator has modulus of continuity of order $h\left(1+\log\frac{D}{h}\right)$, while the remainder operator is Lipschitz continuous. These estimates imply Dini continuity of the scalar flux, providing the regularity needed to control the singular integral arising from differentiation of the transport potential.

The Dini estimate allows the differentiated principal potential to be separated into a bounded interior fluctuation term and an explicit boundary contribution. The latter generates the leading logarithmic singularity, whereas the differentiated remainder contributes only bounded corrections. For every fixed interior angle $0 < \theta < \pi/2$, the resulting corner asymptotics are

$$\partial_{x_k}\phi(re_\theta) = A\log\frac{1}{r} + O_\theta(1), \quad k = 1, 2,$$

and

$$\phi(re_\theta) - \phi(P) = A(\cos\theta + \sin\theta) r \log\frac{1}{r} + O_\theta(r),$$

where $A$ depends on both the scalar flux and the external-source value at the corner.

Thus, the inverse-distance component of the transport kernel determines the leading singular behavior, with the coefficient containing contributions from both scattering and the external source.

The numerical results in Section 6 provide a finite-$r$ consistency test of the predicted $r\log(1/r)$ scalar-flux correction and its analytical coefficient. They also illustrate the influence of subleading terms and discretization errors on finite-$r$ estimates of the leading coefficient.

The principal conclusion is that a rectangular corner preserves the logarithmic first-derivative singularity already present in smooth-boundary transport theory rather than generating a stronger singularity: the leading singularity remains logarithmic. The analysis further identifies its coefficient and the corresponding local scalar-flux expansion. Together, the Dini regularity estimate and fixed-angle asymptotics provide a rigorous basis for characterizing corner effects in two-dimensional radiation transport and their influence on numerical approximation.

**References**

1. D. Wang, "On Solution Regularity of the Two-Dimensional Radiation Transfer Equation and Its Implication on Numerical Convergence," *Journal of Computational and Theoretical Transport*, vol. 52, no. 6, pp. 452–474, 2023. https://doi.org/10.1080/23324309.2023.2292040.

2. E. E. Lewis and W. F. Miller, Jr., *Computational Methods of Neutron Transport*, New York, NY: Wiley, 1984; 2nd ed., La Grange Park, IL: American Nuclear Society, 1993.
3. H. G. Kaper, “A Collection of Problems in Transport Theory,” *Transport Theory and Statistical Physics*, vol. 4, no. 1, pp. 125–134, 1975. https://doi.org/10.1080/00411457508214535.
4. V. A. Kondrat'ev, “Boundary Value Problems for Elliptic Equations in Domains with Conical or Angular Points,” *Transactions of the Moscow Mathematical Society*, vol. 16, pp. 227–313, 1967.
5. P. Grisvard, *Elliptic Problems in Nonsmooth Domains*. Philadelphia: Society for Industrial and Applied Mathematics, 2011.
6. D. Wang, “On Kaper’s Corner Singularity Problem for the Two-Dimensional Radiation Transport Equation,” *Proc. M&C 2027*, Columbus, OH, USA, April 18–22, 2027.
7. K. M. Case and P. F. Zweifel, *Linear Transport Theory*, Reading, MA: Addison-Wesley, 1967.
8. F. W. J. Olver et al., *NIST Digital Library of Mathematical Functions*, National Institute of Standards and Technology, available at https://dlmf.nist.gov.
9. S. G. Mikhlin, *Multidimensional Singular Integrals and Integral Equations*. Oxford, UK: Pergamon Press, 1965.
10. E. M. Stein, *Singular Integrals and Differentiability Properties of Functions*. Princeton, NJ: Princeton University Press, 1970.